\documentclass[10pt,bezier]{article}
\usepackage{amsmath}
\usepackage{amsfonts,amsthm,amssymb}
\usepackage{amsfonts}

\usepackage{amsmath}
\usepackage{amsfonts,amsthm,amssymb}
\usepackage{amsfonts}
\usepackage{enumerate}
\usepackage{amscd}
\usepackage{amsmath,amsfonts,amssymb,amscd}
\usepackage{indentfirst,graphicx,epsfig}
\usepackage{graphicx,psfrag}
\input{epsf}

\begin{document}
\title{Algorithm on rainbow connection for maximal outerplanar graphs\thanks{
The project supported partially by CNNSF (No.11401181)and the first author is supported by Tian Jin Normal University Project (No.52XB1206).}}
\author{Xingchao Deng$^a,$  Hengzhe Li$^b,$ Guiying Yan$^c$ \\
 \footnotesize $^a$College of Mathematics, Tian Jin Normal University\\
 \footnotesize Tianjin City, 300071, P. R. China\\
 \footnotesize  E-mail: dengyuqiu1980@126.com(xcdeng@mail.tjnu.edu.cn )\\
   \footnotesize $^b$College of Mathematics and Information Science, Henan Normal University\\
       \footnotesize Xinxiang 453007, P. R. China\\
       \footnotesize $^c$Academy of Mathematics and Systems Science, CAS, \\
       \footnotesize  Beijing 100190, P. R. China\\}
\date{}
\maketitle
\begin{abstract}
 In this paper, we consider rainbow connection number of maximal outerplanar graphs(MOPs) on algorithmic aspect. For the (MOP) $G$, we give sufficient conditions to guarantee  that $rc(G) = diam(G).$ Moreover, we produce the graph with given diameter $d$ and give their rainbow coloring in linear time. X.Deng et al. $\cite{XD}$ give a polynomial time algorithm to compute the rainbow connection number of MOPs by the Maximal fan partition method, but only obtain a compact upper bound. J. Lauri $\cite{JL}$ proved that, for chordal outerplanar graphs given an edge-coloring,  to verify whether it is rainbow connected is NP-complete under the coloring, it is so for MOPs. Therefore we construct Central-cut-spine of MOP $G,$ by which we design an algorithm to give a rainbow edge coloring with at most $2rad(G)+2+c,0\leq c\leq rad(G)-2$ colors in polynomial time.\\

\end{abstract}

\section{Introduction}

    Graphs considered are finite, simple and connected in this paper.
Notations and terminologies not defined here, see West
\cite{WT}. The concept of rainbow concetion was introduced by Chartrand, Johns, McKeon and Zhang in
2008 $\cite{CJMZ}$. Let $G$ be a nontrivial finite simple connected graph on
which is assigned a coloring $c: E(G)\rightarrow \{1,2,\cdots,n\},
\ n\in \mathbb{N},$ where adjacent edges may have same color. A {\emph rainbow path} in $G$ is a path
with different colors on it. If for any two vertices
of $G,$ there is a rainbow path connecting them, then $G$ is called
{\emph rainbow connected} and $c$ is called a
{\emph rainbow coloring}. Obviously, any $G$ has a trivial rainbow
coloring by coloring each edge with different colors. Chartrand et al.\cite{CJMZ} defined the {\emph rainbow connection number} $rc(G)$ of graph $G$
as the smallest number of colors needed to make $G$ rainbow connected.
For any two vertices $u$ and $v$ in $G,$ the length of a shortest path between them is their distance,
denoted by $d(u, v).$  The eccentricity of a vertex $v$ is $ecc(v) :=max_{x\in V (G)} d(v, x).$
The diameter of $G$ is $diam(G) := max_{x\in V (G)} ecc(x).$
The radius of G is $rad(G) := min_{x\in V (G)} ecc(x).$ Distance between a vertex $v$ and a set
$S \subseteq V (G)$ is $d(v, S) := min_{x\in S} d(v, x).$ The k-step open neighbourhood of a set $S \subseteq V (G)$
is $N_k(S):= \{ x \in V (G)|d(x, S) = k\}, k \in \{0, 1, 2, \ldots \}.$ The degree of a vertex $v$ is
$degree(v):= \mid N_1({v})\mid.$ The maximum degree of $G$ is $\Delta(G) := max_{x\in V (G)} degree(x).$ The girth of a graph $G$ is $g(G):=\mbox{the length of maximal induced cycle in $G$}.$ A vertex is called pendant if its degree is $1.$ Let $n(G)=\mid V(G)\mid$ and $e(G)$ be the size of $G.$ Obviously $diam(G)\leq rc(G) \leq e(G).$ From \cite{CJMZ}, we know that rainbow connection number of
any complete graph is 1 and that of a tree is its size.

Obviously, we know that cut-edges must have distinct colours when G is rainbow connected. Thus stars have arbitrarily large rainbow connection number while having diameter 2. Therefore, it is significant to seek upper bound on $rc(G)$ in terms of $diam(G)$ in 2-edge-connected graphs. Chandran et al.\cite{CDRV} showed that $rc(G)\leq rad(G)(rad(G)+2)$ when G is 2-edge-connected, and hence $rc(G)\leq diam(G)(diam(G)+2).$ Li et al. \cite{LLL} proved that $rc(G)\leq 5$ when G is a 2-edge-connected graph with diameter 2. Li et al. \cite{LLS} proved that $rc(G)\leq 9$ when $G$ is a 2-edge-connected graph with diameter 3.

Recalling an {\emph outerplanar graph} is a planar graph which has a plane embedding
with all vertices placed on the boundary of a face, usually taken to
be the exterior one. A MOP is an outerplanar graph which can not be added any line without losing
outerplanarity.

By \cite{LB}, a MOP can be recursively defined as follows:
$\mathbf{(a)}$ $K_3$ is a MOP.
$\mathbf{(b)}$ For a MOP $H_1$ embedded in the plane
    with vertices lying in the exterior face $F_1,$  $H_2$ is obtained by joining a new vertex to two adjacent vertices on $F_1.$  Then $H_2$ is a MOP.
$\mathbf{(c)}$ Any MOP can be constructed by finite steps of $\mathbf{(a)}$ and $\mathbf{(b)}.$

Note each inner face of a MOP $H$ is a triangle and the connectivity $\kappa (H)=2.$ Moreover, $H$ can be represented by two line arrays
$High(1),High(2),\cdots,High(n)$ and
$Low(1),Low(2),\cdots,Low(n).$ Here for any vertex $i,$ $High(i)$ and $Low(i)$ are labels of its two neighbors whose labels are
less than $i,$ and $High(i)> Low(i);$ and $High(1),
Low(1)$ and $Low(2)$ are undefined, and $High(2)=1.$
Figure 1 illustrates a MOP and its canonical representation.

\begin{figure}[ht]
\begin{center}
\includegraphics[width=6cm,totalheight=30mm]{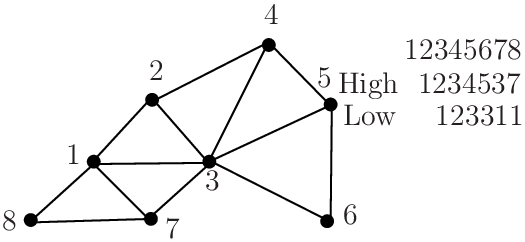}
\end{center}
\caption{$\mathrm{Example}$}\label{fig1}
 \end{figure}


\noindent{\bf Property (A) } A graph is outerplanar if and only if it has no $K_4$ or $K_{2,3}$ minor.\\

 We summarize some results for the rainbow connection number of graphs in the following.

Huang et al. proved that if $G$ is a bridgeless outerplanar graph of order $n$ and $diam(G) = 2,$ then $rc(G) \leq 3$ and the bound is tight. Moreover they proved that if $diam(G) = 3,$ then $rc(G) \leq 6,$ in $\cite{XH}$.

\vspace{1mm}

\noindent\textbf{Theorem $1.1^{\cite{CJMZ}}.$ } For cycle $C_n$, we have
\begin{equation*}
rc(C_n)=\left\{
\begin{array}{ll}
 \frac{n}{2} & {\rm if~}n\ {\rm is\ even,}\\
\lceil\frac{n}{2}\rceil & {\rm if}~n~{\rm is\ odd.}
\end{array}
\right.
\end{equation*}%
  Chandran et al. \cite{CDRV} studied the
relation between rainbow connection numbers and connected dominating
sets, and they obtained the following results:\\

$\mathbf{(1)}$ For any bridgeless chordal graph $G,$ $rc(G)\leq 3 rad(G).$ Moreover, the result is tight.

$\mathbf{(2)}$ For any unite interval graph $G$ with $\delta(G)\geq 2,$ $rc(G)=diam(G).$\\

A finite simple connected graph $G$ is called a {\it Fan} if it is $P_n\vee K_1$ (the join of $P_n$ and $K_1$), denoted by $Fan_{n},$ for some $n\in\mathbb{N}\setminus\{1\}.$ Here the vertex $v$ of $K_1$ is
called {\emph central vertex,} the edges $v_iv_{i+1}(1\leq i\leq n-1)$ of
$P_n=(v_1v_2\cdots v_n)$ are called {\emph path edges,} and the edges
$v_iv$ between $P_n$ and $K_1$ are called {\emph spoke edges.} \\

\noindent{\bf Theorem $1.2^{\cite{XD}}.$} The rainbow connection number of $Fan_{n}$ satisfies

\begin{equation*}
rc(Fan_n)=\left\{
\begin{array}{ll}
1 & {\rm if}~~n=2,\\
2 & {\rm if}~~3\leq n\leq 6,\\
3 & {\rm if}~~n\geq 7.
\end{array}
\right.
\end{equation*}

\vskip 2mm
%
%
%
%
%
%

\noindent{\bf Theorem $1.3^{\cite{XD1}}.$} Let G be a bridgeless outerplanar graph of order n.\\

   1.  If $diam(G)=2,$
   then \begin{equation*}
rc(G)=\left\{
\begin{array}{ll}
3 & {\rm if}~~G=F_{n}~(n\geq 7)~ {\rm or~ C_5,}\\
~~\\
2 & {\rm otherwise}.
\end{array}
\right.
\end{equation*}

   2. If $diam(G) = 3,$ then $3\leq rc(G) \leq 4$ and the bound is tight.\\

Following, we give a theorem on edge, vertex cut set and their rainbow connection for a connected graph $G.$\\

\noindent{\bf Theorem $1.4.$} Let $G$ be a connected graph and $S_1,$ $S_2 $ be two disjoint edge cuts, then they must
be colored by at least two different colors in order to make $G$ rainbow connected.\\

\noindent{\bf Proof.} Let $V_{i1},V_{i2}$ be the end vertex sets of $S_i$ $(i=1,2),$ then  $V_{i1},V_{i2}$ are
vertex cuts corresponding to $S_i,$ $(i=1,2)$ and $X_i,$ $\overline{X_i}$ be the vertex sets separated
by $S_i,$ $(i=1,2).$ Clearly, there are no red edges in the graph $G$ showing as in Figure 2.
Since $S_2$ is an edge cut of $G,$ no edges between $V_{11}\bigcap (X_2\setminus V_{21})$ and $V_{12}\bigcap (\overline{X_2})$ and so no edges between $V_{12}\bigcap (X_2\setminus V_{21})$ and $V_{11}\bigcap (\overline{X_2}) ,$  between $V_{11}\bigcap (\overline{X_2}\setminus V_{22}),$  $V_{12}\bigcap X_2 ,$   between $V_{12}\bigcap (\overline{X_2}\setminus V_{22})$ and $V_{11}\bigcap X_2 . $

If $(\overline{X_1}\setminus V_{12} )\bigcap (\overline{X_2}\setminus V_{22})$ and $(X_1\setminus V_{11}) \bigcap (X_2\setminus V_{21})$ are not empty, thus  must
 have rainbow path through $S_1$ and $S_2$ in order to make $G$ rainbow connected, therefore Theorem 1.4 is correct. Other cases can be proved by the same method.                             $\hfill \Box$\\

\begin{figure}[ht]
\begin{center}
\includegraphics[width=3cm,totalheight=30mm]{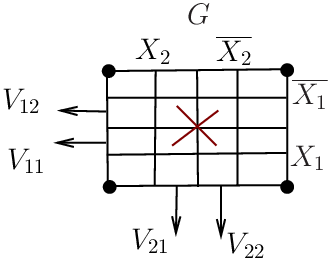}
\end{center}
\caption{Partition of $V$ determined by the edge cuts $S_1$ and $S_2$ }\label{fig1.4}
\end{figure}

\vskip 2mm

\section{rc(G) of minimum maximal outerplanar graph(MMOP) G}

 We give a sufficient condition for some graphs of maximal outerplanar graphs whose rainbow connection numbers equal their strong rainbow
connection numbers, just their diameters in this section.

The MMOP $G$ is a MOP with given diameter having minimum vertices. Recall a MOP is outerplanar, therefore $G$ does not contain $K_4$ or
$K_{2,3}$ minor (\cite{GC}).  We construct the MMOPs with diameter $2$ by its recursive definition through the following three steps.\\

\noindent{\bf Step $1.$} $K_3$ with vertices $a,b,c$(as Figure 3 showing) as the start of the process constructing the MMOPs with diameter $2.$

\noindent{\bf Step $2.$} By the symmetry  of $K_3$, we could add one new vertex $d$ joining any two vertexes of $a,b,c $ (For example as Figure 4 showing).
\begin{figure}[ht]
\begin{center}
\includegraphics[width=3cm,totalheight=20mm]{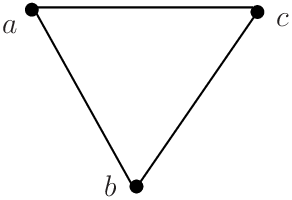}
\end{center}
\caption{$\mathrm{H_1=K_3}$}\label{fig2}
\end{figure}



\noindent{\bf Step $3.$} Since  Figure 4 has diameter $2$ and only $a,d$ has distance two, in the process of the construction we add a new vertex $e$ and joining any one of edges of exterior face of  Figure 4. By symmetry, we obtain the following  MOP(Figure 5) with diameter $2.$
\begin{figure}[ht]
\begin{center}
\includegraphics[width=3cm,totalheight=20mm]{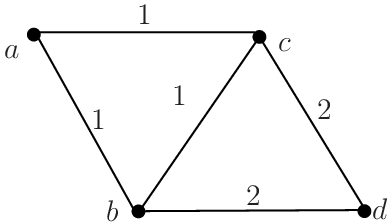}
\end{center}
\caption{$\mathrm{K_3+d}$}\label{fig3}
\end{figure}

\begin{figure}[ht]
\begin{center}
\includegraphics[width=3cm,totalheight=20mm]{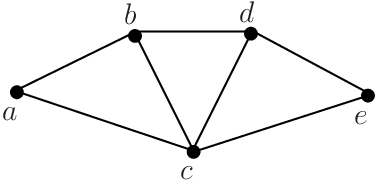}
\end{center}
\caption{$\mathrm{F_4}$}\label{fig4}
 \end{figure}

 Since Figure 4 and Figure 5 have diameter 2 and the resulting graph of deleting any one of the
  vertices which have distance 2 in Figure 4 is $K_3$ with diameter $1,$ thus Figure 4 is the  MMOP with diameter $2.$

 In Figure 5, the two vertex pairs $a,d$ and $a,e$ have distance $2.$ We add a new vertex adjacent to $d$ and $e$ and obtain an outerplanar graph with diameter $3.$\\

\begin{figure}[ht]
\begin{center}
\includegraphics[width=3cm,totalheight=20mm]{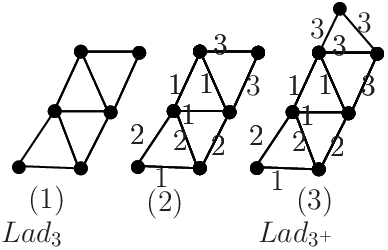}
\end{center}
\caption{$\mbox{diameter-3~MMOPs}$}\label{fig5}
 \end{figure}

Note that $Lad_3$ and $Lad_{3^+}$ have diameter 3 and $Lad_3$ is the MMOP with diameter $3.$ Following, we call the MMOP with diameter $d$ $Lad_d,$  $Lad_{d^+}$ is obtained by adding a new vertex to $Lad_d$ as $Lad_3$ and $Lad_{3^+} $ showing.\\


\noindent{\bf Theorem $2.1.$} If $G$ is $Lad_d$ or $Lad_{d^+}$ with diameter $d,$ then $rc(G)=src(G)=d.$ \\

\noindent{\bf Proof.} We prove the theorem by induction method.
Through the above three steps constructing the MMOPs with diameter 2 and 3, we know that the theorem is corrected when $d=2,3.$

Now suppose that we obtain the MMOP $G$ with diameter $d-1$ and $rc(G)=src(G)=d-1.$
 Then $G$ has a vertex pair, denoted by $a,b,$  having distance $d-1$ and there is a $d-1-path$ connecting them.
 Since $G$ is a MMOP with diameter $d-1,$ therefore there is only one vertex pair with distance $d-1.$

 Following, we construct the MMOP with diameter $d.$ $G$ has an $ab-path$ with distance $d-1$
 and is a  MOP, then we add a new vertex $e$ to $G$ adjacent to $b$ and it's preceding vertex $c,$ on the $ab-path.$ The edges $ec,eb$ are colored by $d-1$ respectively. Finally, we add another vertex $f$ to the above graph and adjacent it to $b,e.$  The new graph has diameter $d$ and has minimum vertex. The edges $fe,fb$ are colored by $d.$ We give the construction process in Figure 7.

 Now, we prove the colouring of Figure 7 is a rainbow and strong rainbow coloring of the graph. The part $A$ of Figure 7 is a MMOP with diameter $d-1$ and has a strong rainbow coloring with $d-1$ colors. When $e$ is added to $A,$  we connected it
 to $b,c$ which are colored by $d-1$ with the same color of the edge $cb.$  Therefore there is rainbow path between $e$ and other vertices of $A$ through the edge $ce.$ Obviously, the two vertex pairs $ec,eb$ are strong rainbow connected by the coloring. By the same reason, we can prove the resulting graph added the vertex $f$ and $fe,fb$ colored by $d$ is strong rainbow connected, which is $Lad_d.$ Moreover, we add a new vertex $g$ which is adjacent to $e$ and $f,$  the edges $ge,$ $gf$ is colored by $d.$ The resulting graph is $Lad_{d^+}.$  They can be obtained by the following algorithm.      $\hfill \Box$\\

\begin{figure}[ht]
\begin{center}
\includegraphics[width=3cm,totalheight=20mm]{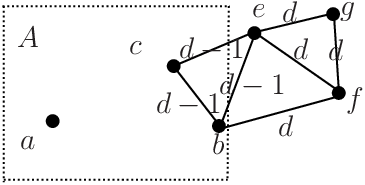}
\end{center}
\caption{$\mbox{Producing ~graph~ with ~diameter ~d}$}\label{fig5}
 \end{figure}

\newpage
\noindent{\bf Algorithm 1: Giving the graphs($Lad_d,$ $Lad_{d^+}$) with diameter $d$ and their rainbow(strong rainbow) colorings}\\

\noindent{\bf Input:} A number $d\geq 4.$\\
\noindent{\bf Output:} The MMOP $G$ with diameter $d$ and a strong rainbow coloring of it.

{\bf Step 1:} Given a $K_3$ with vertices $1,2,3,$ $d(1,1)=0,d(1,2)=1,d(1,3)=1.$ Adding vertex $4$ to $K_3$

\hskip0.2cm and $4$ is adjacent to $2,3,$ then $d(1,4)=2.$ Now, adding a new vertex $5$ to the above graph with

\hskip0.2cm vertices $\{1,2,3,4\},$  which is adjacent to $4$ and it's preceding neighbor on a $2-length-path,$

\hskip0.2cm then $d(1,5)=2.$  Adding a vertex $6$ adjacent to vertices $4,5$ to the graph with vertices $\{1,2,3,4,5\},$

\hskip0.2cm then $d(1,6)=3.$

{\bf Step 2:}  In this step, we construct the MMOP $Lad_{d}$  with diameter $d$ and it's strong rainbow coloring.

\hskip0.2cm{\bf begin}

\hskip0.2cm{\bf for}   $1\leq i\leq 2d$ in $G,$ {\bf if} $d(1,i-1)=d(1,i)-1<d,$ {\bf  add} a new vertex $(i+1)$ connecting to $i-1,i$

\hskip0.2cm {\bf then} $d(1,i+1)=d(1,i),$  {\bf Color the edges} $(i-1,i+1)(i,i+1)~by~i+1;$

\hskip0.2cm {\bf do} $i\leftarrow i+1.$

\hskip0.2cm {\bf if} $d(1,i-1)=d(1,i)<d,$ {\bf  add} a new vertex $(i+1)$ connecting to $i-1,i$ {\bf then}

\hskip0.2cm$ d(1,i+1)=d(1,i)+1,$  {\bf Color the edges} $(i-1,i+1)(i,i+1)~by~i-1;$

\hskip0.2cm {\bf do} $i\leftarrow i+1.$

{\bf Step 3:} Adding a vertex $2d+1$ to the above graph $G$ with $2d$ vertices and connecting it to $2d-1,2d.$

\hskip0.2cm Coloring the new edges $(2d+1,2d),(2d+1,2d-1)$ by $2d,$ then we obtain the graph $Lad_{d^+}$ and it's

\hskip0.2cm strong rainbow coloring.

\hskip0.2cm {\bf end}\\

\noindent{\bf Problem $ 2.2^{\cite{XL}}.$} Characterize those graphs $G$ with $rc(G) = diam(G),$ or give some
sufficient conditions to guarantee that $rc(G) = diam(G).$ Similar problems for the
parameter $src(G)$  can be proposed.\\

Obviously, it is very difficult to give the sufficient conditions for general graphs. Up to now, we only know very few graphs whose rainbow connection numbers and their strong rainbow connection numbers are their diameters, which are particular graph classes. The MMOP is a graph class whose rainbow and strong rainbow connection numbers are their diameters. Figure 8 gives three  MOPs with diameter $3$ and rainbow coloring using $3$ colors. Therefore the condition of {\bf theorem 2.1} is not necessary.

\begin{figure}[ht]
\begin{center}
\includegraphics[width=5cm,totalheight=20mm]{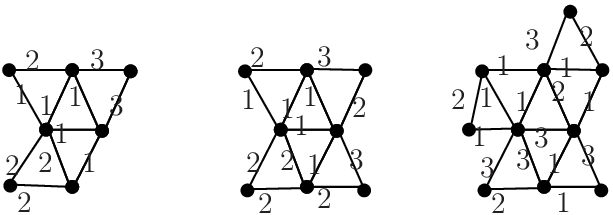}
\end{center}
\caption{$\mbox{diameter~3~maximal~outerplanar}$}\label{fig6}
 \end{figure}

\section{Central-cut-spine of MOP}

In this section, Our main result is to give the definition of {\emph Central- cut-spine} of  MOP, which play a key role in algorithm 3. Following we give algorithm 2 to compute the Central-cut-spine of a MOP.

 A graph $G$ is called {\emph 2-degenerate} if any of its subgraph has a vertex with degree 2 or less.


\noindent{\bf Theorem $3.1^{\cite{GH}}.$} For a connected
graph $G$ with at least 2 vertices, it is an MOP iff
the following hold: (a) for any vertex $v$ of $G,$ its
neighbors induce a path in $G,$ (b) $G$ is 2-degenerate.\\

It is well known that a MOP $G$ can be embedded in the plane such that every vertex lines on the boundary of the exterior face, all exterior edges form a
{\emph Hamiltonian cycle} $C=(v_1v_2\cdots v_nv_1).$ In $G ,$ a {\emph Hamiltonian degree sequence} $D=(d_1,d_2,\cdots, d_n,d_1)$ is the degree sequence of vertices $v_1,v_2,\cdots,v_n,v_1.$
In $\cite{TB},$ T. Beyer et al.  gave an {\emph Algorithm} which takes a {\emph Hamiltonian degree sequence} and produces the unique corresponding MOP in linear time. For any edge $v_{s}v_{t}$ of $G,$ which is a chordal edge of $G$ when $\mid s-t\mid\geq2,$ the two vertices are a {\emph 2-vertex-cut} set. Since $g(G)=3,$ any two vertices incident to a chord of $G$ have exactly two common neighbors, while two vertices incident to an outer edge of $G$ have exactly one common neighbor. Now if $G$ has order $n.$ Using the fact that the boundary of the exterior region of $G$ is a hamiltonian cycle and the boundary of every interior region of $G$ is
a triangle, which follows that $G$ has $2n - 3$ edges and $n - 1$ regions by Euler's formula. Thus $G$ has
$n - 3$ chords and $n - 2$ interior triangles.\\

\noindent{\bf Theorem $3.2^{\cite{TB}}.$} A MOP $G$ is determined uniquely up to isomorphisms by its Hamiltonian degree sequence $D=(d_1,d_2,\cdots, d_n,d_1).$\\

A graph is chordal if every cycle of length greater
than three has a chord, which is meaning that there is an edge joining two nonconsecutive vertices of the cycle.\\

\noindent{\bf Theorem $3.3^{\cite{CG}}.$}  $2rad(G)-2 \leq diam(G) \leq 2rad(G)$ for any connected chordal graph $G.$ Moreover, if $2rad(G) - 2 = diam(G),$ then $G$ has a 3-sun as an induced subgraph.\\

Let $\eta(G)$ be the smallest integer such that every edge of $G$ belongs to a cycle of length at most $\eta(G).$ \\

\noindent{\bf Theorem $3.4^{\cite{XHH}}.$} For every bridgeless graph $G,$
$rc(G)\leq \sum^{rad(G)}_{i=1}min\{2i + 1, \eta(G)\} \leq rad(G)\eta(G).$\\

The Central-cut-spine(CCS) of a MOP $G$ is a tree generated by the following
method: First, we choose an center vertex, denoted by $v,$ marked red, on the unique hamiltonian cycle of $G.$
Second, give the sets  $N_i(v):= \{ x \in V (G)|d(x, v) = i\}, i \in \{0, 1, 2, \ldots ,rad(G) \}.$
Third, for $1\leq i\leq rad(G)-1,$ if $\mid N_i(v)\mid=2 $ and they are adjacent, then shrink the corresponding
edge to obtain a new vertex, marked green, adjacent to the new vertex obtained by the vertices in $N_{i-1}(v)$ pertinent to $N_i(v); $ If  $\mid N_i(v)\mid\geq 3,$  then we partition them into $N_{ij}(v), 1\leq j \leq \mid N_i(v)\mid ,$ for each $j,$ the vertices of $N_{ij}(v)$ forming all the path edges of a Fan-structure, for any edge $v_{i_s}v_{i_t},$  $\mid i_s- i_t\mid\geq2$ and
 they have neighbors in $N_{i+1}(v),$ then shrink the corresponding edge to obtain a new vertex, named $v_{ijv_{i_s}v_{i_t}}$ marked green, and adjacent
to the new vertex $v_{i-1j^{'}v_{(i-1)_1}v_{(i-1)_2}}$ pertinent to $v_{i_s}v_{i_t}.$\\



\begin{figure}[ht]
\begin{center}
\includegraphics[width=6cm,totalheight=30mm]{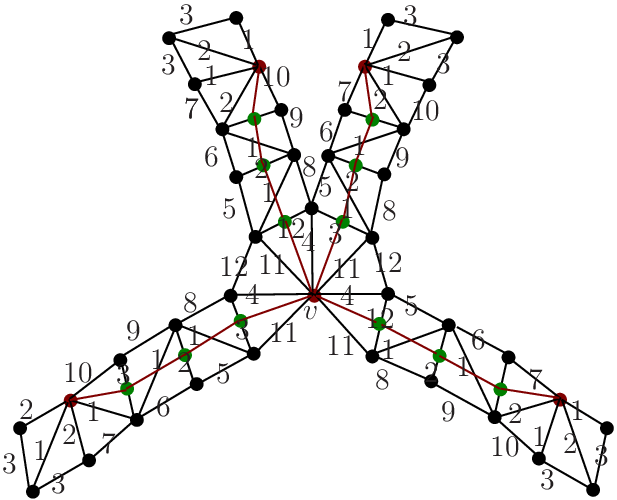}
\end{center}
\caption{$\mathrm{Fan_7+4-sun+Lad_4}$}\label{Fan70}
 \end{figure}

 An example:  we discuss the property of Figure \ref{Fan70} and give it's CCS. Notice that it has $40$ vertices and the center vertex $v,$ radius $5$ and diameter $10.$ By Theorem 3.3, we know that  $rc(Figure~\ref{Fan70})\leq 15.$ In Figure \ref{Fan70}, we give it a $12-rainbow ~coloring.$

 The red vertices, edges and green vertices constitute it's CCS which is computed by the following operations: We choose an eccentricity vertex, denoted by $v,$ with maximum degree and  the vertices $N_1(v).$  Choose the adjacent vertices in $N_1(v) $, which are cut sets of $G$ and replaced by new vertices, which are marked green and adjacent to $v.$  We proceed the same operation on $N_i(v),2\leq i\leq 4$ as $N_1(v).$  The last step we choose the vertices, in  $N_4(v)=N_{radius-1}(v),$  whose neighbors forming a Fan-structure with them as central vertices, which are marked red and adjacent to the preceding new vertices. Now we obtained the Central-cut-spine of Figure \ref{Fan70}, denoted by CCS$(Figure\ref{Fan70}).$

\begin{figure}[ht]
\begin{center}
\includegraphics[width=6cm,totalheight=40mm]{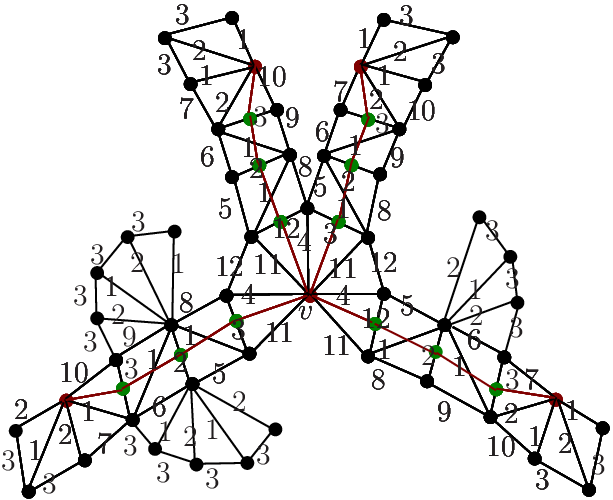}
\end{center}
\caption{$\mathrm{Fan_7+4-sun+FLad_4}$}\label{FFan7}
 \end{figure}

  Furthermore, we consider the $Lad_4s$ and their edge cuts of Figure \ref{Fan70}. By theorem 1.4, in order to give Figure \ref{Fan70} a rainbow coloring, we know that every edge cut has a color which is different from other edge-cuts having in $Lad_4.$ Since it has four $Lad_4s,$ at least two $Lad_4s$ are colored by eight different colors. Obviously, $rc(Figure\ref{Fan70})\geq 10.$ When the $Lad_4s$ are replaced by $Lad_ns,$ the resulting graph, denoted by $Fan_7+4-sun+Lad_n,$ has diameter $2n+2$ and it's rainbow connection number is at least $2n+2.$  With the same method of coloring in Figure \ref{Fan70}, we can rainbow color $Fan_7+4-sun+Lad_n$ by $2(n-1)+6=2n+4,$ which is {\bf less $n-1$} than the $3rad(Fan_7+4-sun+Lad_n).$ As $n\rightarrow \infty,$ we know that the upper bound $3rad$ in Theorem 3.3 may be arbitrarily far from $2n+4.$ Therefore, the upper bound is not good for rainbow connection number of MOPs. Figure \ref{FFan7} gives an another graph with the Central-cut-spine formed by red edges, vertices and green vertices, which has a rainbow coloring with 12 colors.

 \begin{figure}[ht]
\begin{center}
\includegraphics[width=5cm,totalheight=40mm]{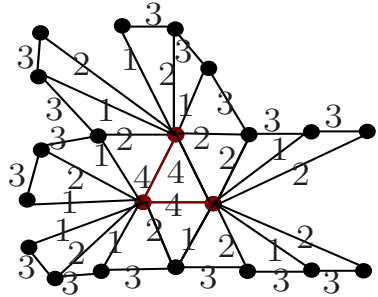}
\end{center}
\caption{$\mathrm{diameter-3+a9b9c9}$}\label{fig8}
 \end{figure}

 %

Figure \ref{fig8} gives an example of MOP with diameter $3,$ where the maximal Fan-structures have $10$ vertices and  the red vertices and edges give it's CCS$(Figure\ref{fig8}).$ The Fan-structures can be replaced by any $F_n$ with $n\geq 2,$ and the resulting graph $G$ can be rainbow colored by at most $4$
colors with the same method as Figure \ref{fig8} showing. As the above description, we know that the MOPs having any length hamiltonian cycles. We may color the hamiltonian cycle making it rainbow connected by $\frac{1}{2}n(G)$ colors. However, as $n\rightarrow \infty,$ the number of colors $\rightarrow \infty.$ Thus, rainbow coloring the hamiltonian cycle is not an effective method for MOPs. Moreover, when we compute the upper bound of $rc(G)$ by the method of Maximal fan partition method showed in $\cite{XD},$  the colors needed will increase to 9 as $n$ increasing. The method is not effective either.\\

 Following, we give an effective algorithm to compute an upper bound of rainbow connection number for  MOP with the help of it's CCS.\\

\noindent{\bf Theorem 3.5.} If $G$ is a MOP and $v_r$ is the root of CCS$(G),$ then we have the
following results:
 \begin{enumerate}
   \item {Each $v_r-path$ in CCS$(G)$ corresponding to two edge disjoint paths in $G.$}
   \item {The green vertices of CCS$(G)$ are corresponding to cut edges or the edges connecting to the vertices which are central vertices in the construction of CCS$(G).$}
   \item {If $v_i$ is a red vertex of CCS$(G),$ then it is a vertex of $G$ and a central vertex
   in the construction of CCS$(G).$}
 \end{enumerate}

\noindent{\bf Proof:} 1. By the construction of the CCS$(G)$ of $G,$ we know that $v_r$ is a red vertex and the root of CCS$(G).$

Case1: The leaf vertex is a red vertex in CCS$(G),$ then it is a vertex of $G.$ Since $G$ is a $2-connected$ graph, there are two edge disjoint paths between it and $v_r$ in $G,$ one of which has length at most $rad(G)-1.$  Since all edges on the path are situated in triangles, there is another path  with length no more than $2(rad(G)-1).$

Case2: If the leaf is a green vertex in CCS$(G),$ then it is corresponding to two adjacent vertices which are a vertex cut set and situated in  $N_{(rad(G)-1)}$ or $N_{(rad(G)-2)}.$ Therefore there are paths, with length at most $rad(G)-1,$ between $v_r$ and the  two adjacent vertices which are corresponding to the leaf.  If they are edge disjoint then we choose them as the corresponding two
paths, otherwise we choose one path with length $rad(G)-1$ or $rad(G)-2.$ Since any edge of the path located in one triangle and every common edge can be replaced by two adjacent edges in the corresponding triangle, thus we can obtain a path between $v_r$ and the another vertex with length at most $(rad(G)-1)+c_1,$ where $c_1$ is the number of common edges on the two paths.

By the construction of CCS$(G),$ we produce the corresponding two paths in $G$ for a $v_r-path$ of
CCS$(G):$
\begin{enumerate}
\renewcommand{\labelenumi}{(\arabic{enumi})}
  \item {$v_r$ is the start point of the two paths;}
  \item {If the successor of $v_r$ is a green vertex on one $v_r-path,$ then whose corresponding vertices in $G$ are the successors of the two paths respectively; }
  \item {Now we consider the vertex having distance two to $v_r.$ If it is a green vertex $v_{{i2}v_{i_1}v_{i_2}}$ of CCS$(G),$ which is corresponding to $v_{i_1}$ and $v_{i_2}$ in $G.$ We choose the vertex having maximum degree, w.l.g. $v_{i_1}.$ Since $v_{i_1}$ is situated in $N_2(v_r),$ there is an path between it and $v_r$ with length $2.$ Moreover any edge of the path located in a triangle, we could choose an edge disjoint path with the $2-length$ path between $v_{i_2}$ and $v_r$ with length
      at most $4.$  Following we produces the above process.}
\end{enumerate}

  By the construction of CCS$(G),$ 2. 3. are correct Obviously.   $\hfill \Box$\\


A. Farley et al. introduced the notion of edge eccentricities by the separation property of an edge for outerplanar graphs in \cite{FP}.\\

\begin{figure}[ht]
\begin{center}
\includegraphics[width=5cm,totalheight=30mm]{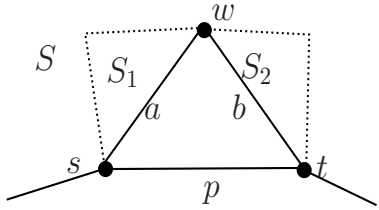}
\end{center}
\caption{$\mathrm{An~ edge~ p = (s, t)~ with~ a~ nonempty~ side~ S.}$}\label{fig9}
 \end{figure}

 \noindent{\bf Definition $3.6^{\cite{FP}}.$} Let $p = (s, t)$ be any edge of an outerplanar graph $G,$  showing as in Figure \ref{fig9}. A. Farley et al. defined four values $e(p, x, S),$ which are called edge eccentricity of $p,$ one for each vertex $x$ $(s ~\mbox{or}~ t)$ and side $S$ of the edge $p.$ The absolute value of $e(p, x, S)$ equals the eccentricity of the vertex $x$ in
the induced subgraph $G\left[S \cup \{s, t\}]\right.$  The value $e(p, x, S)$ is negative iff all
vertices of $S \cup\{s, t\}$ at distance $d = e(p, x, S)$  from $x$ lie at distance $d - 1$ from the other end vertex of $p.$\\

In \cite{FP}, A. Farley et al. gave an {\emph edge eccentricity algorithm:
}
  Given a MOP $G,$ calculate the eccentricities of its edges as follows:\\

(a) For all edges $p = (s, t) $ on the Hamiltonian cycle of $G,$ assign the value $- 1$
to $e (p, x, \emptyset) ,$ where $x \in  \{s,t\} .$\\

(b) For each triangle $(s, t, w),$ the values $e(a, s, S_1),$ $e(a, w, S_1),$ $e(b, t, S_2),$
and $e(b, w, S_2)$ are defined, one assign values of $e(p, s, S)$ and $e(p, t, S)$ according to
the following rules:\\

Given an edge $p = (s, t)$ of a MOP with a non-empty side $S,$ let $e_1, e_2$
and $r$ represent the values of $e(b, w, S_2),$ $e(b, t, S_2)$ and the eccentricity of $s$ in the
subgraph $G[S_2 \cup\{ s, t, w\}],$ respectively.\\

\begin{equation*}
 r=\left\{
\begin{array}{ll}
  -(1+ e_2)  & \  e_2 > 0, \\
 \mid e_2\mid,  & {\rm otherwise.}
\end{array}
\right.
\end{equation*}

Given an edge $p = (s, t)$ with a non-empty side $S,$ let $e_3$ and $d_1$
represent the values $e(a, s, S_1)$ and $e(p, s, S),$ respectively. Let $r$ be the eccentricity of
$s$ in the graph $G[S_2\cup\{s, t, w\}].$  \\

\begin{equation*}
 d_1=\left\{
\begin{array}{ll}
    e_3  & \  \mid e_3\mid \geq \mid r\mid, \\
 r,  & {\rm otherwise.}
\end{array}
\right.
\end{equation*}

 If an edge $p = (s, t)$ with a non-empty side $S,$\\

 \begin{equation*}
q=\left\{
\begin{array}{ll}
    -(1 + e(a, s, S_1))  & \ {\rm if~} e(a, s, S_1)>0, \\
  \mid e(a, s, S_1)\mid  & {\rm otherwise.}
\end{array}
\right.
\end{equation*}

\begin{equation*}
  {\rm The ~eccentricity~}
d_2 = e(p, t, S)=\left\{
\begin{array}{ll}
\mid e(b, t, S_2)\mid &\ {\rm if~}\mid e(b, t, S_2)\mid\geq \mid q\mid,\\
\mid e(b,t, S_2)\mid & {\rm otherwise.}
\end{array}
\right.
\end{equation*}

By the above algorithm, A. Farley et al. give the following algorithm, denoted by Farley DRC algorithm, to compute the diameter
and center for outerplanar graph $G.$\\

(1) For every vertex $u$ of $G,$ choose an edge $p$ incident with $u,$ then $ecc(u) :=$ the maximum of the absolute values of the two pertinent edge eccentricities
of $p.$\\

(2) $diam(G) := max_{x\in V (G)} ecc(x),$ $rad(G) := min_{x\in V (G)} ecc(x).$ \\

(3) The center of $G$ is the vertex set $C(G):=\{v:ecc(v)=rad(G),v\in V(G)\}.$\\

The time complexity of Farley DRC algorithm for computing all vertex eccentricities in
an outerplanar graphs is $O(n),$ where $n=\mid V(G)\mid.$  \\

A vertex with whose neighbors inducing a clique in $G$ is called {\emph simplicial}.  A {\emph simplicial  elimination ordering} ({\emph perfect elimination ordering}) is a vertex ordering $v_n,\cdots, v_2,v_1$ for which each vertex $v_i$ is {\emph simplicial} in the induced graph by $\{v_1,\cdots,v_i\}.$ It's reverse is called a {\emph simplicial construction ordering} of $G.$ The simplicial construction ordering of a chordal graph $G$ can be found by Maximum Cardinality Search(MCS) in time $O(n(G)+e(G)).$ The MCS algorithm is a simple linear time algorithm
that choose a vertex $x$ and let $f(x) = 1,$ where $f: V(G)\rightarrow
\{1,\cdots,n(G)\}$ is a function, and produces an elimination ordering in reverse. For each vertex $v,$ it maintains an integer weight $l(v)$ that is the cardinality of the already processed neighbors of $v;$ and produces a simplicial construction ordering when a chordal graph is input.\\

Following we give a polynomial time algorithm to compute the Central-cut- spine of MOP.\\
%
%
%

\noindent{\bf Algorithm 2 producing the  Central-cut-spine  of MOP ( CCS-Algorithm):}\\

\noindent{\bf Input:} A MOP $G.$\\
\noindent{\bf Output:} The Central-cut-spine of $G.$

{\bf Step 1:} Finding a  Hamiltonian degree sequence $D=(d_1,d_2,\cdots,d_n,d_1)$ which
is respond to the

vertex sequence $(v_1,v_2,\cdots,v_n,v_1)$ of $G.$

{\bf Step 2:} In this step, we have a simplicial construction ordering of vertices.

\hskip0.2cm{\bf begin}

\hskip0.2cm{\bf for} all vertices $v$ in $G$ {\bf do} $l(v)=0;$

\hskip0.2cm{\bf for} $i = 1$ {\bf up to} $n$ do

\hskip0.2cm Choose an unnumbered vertex $z$ of maximum weight; $f(z) = i;$

\hskip0.2cm{\bf for} all unnumbered vertices $y\in N(z)$
\hskip0.2cm {\bf do} $l(y) = l(y) + 1;$

\hskip0.2cm {\bf end}

{\bf Step 3:} In this step, we have all maximal Fans of $G.$

\hskip0.2cm {\bf begin}

\hskip0.2cm {\bf for} $i=1$ {\bf up to} $n$ {\bf do} $i:$ $N\left(f^{-1}(i)\right)=\left\{f^{-1}(i)\right\};$

\hskip0.2cm{\bf for} $i=1$ {\bf up to} $n$

\hskip0.2cm{\bf if} $\left\vert N\left(f^{-1}(i)\right)\right\vert <d_{f^{-1}(i)}+1$ do

\hskip0.3cm{\bf for} $j=1,\cdots,i-1,i+1,\cdots,n$ do

\hskip0.3cm{\bf if} $f^{-1}(j)$ and $f^{-1}(i)$ are adjacent, $N\left(f^{-1}(i)\right)\leftarrow N\left(f^{-1}(i)\right)\cup
                               \left\{f^{-1}(j)\right\},$ $j\leftarrow j+1;$

\hskip0.3cm otherwise $j\leftarrow j+1;$

\hskip0.2cmelse $i\leftarrow i+1;$

\hskip 0.2cm{\bf end}

{\bf Step 4:} Finding the center vertices, denoted by $C(G),$  of $G$ by the Farley DRC algorithm.

{\bf Step 5:} Choose a vertex in $C(G)$ with minimum degree, denoted by $v_r,$ $1\leq r\leq n.$ Giving a BFS

\hskip0.2cm$v_r-tree$ $T$ in $G$ with predecessor function $p,$ a level function $\ell$ such that $\ell(v) = d_G(v_r, v)$ for all

\hskip0.2cm$v\in V .$

{\bf Step 6:} In this step, we produce the Central-cut-spine of $G.$

\hskip0.2cm {\bf begin}

\hskip0.2cm {\bf Let} $S_{r}$ be the successor function, $P_{r}$ be the predecessor function and $l_{r}$ be level function of CCS$(G).$

\hskip0.2cm $l_{r}(v_r)\leftarrow 0.$

\hskip0.2cm {\bf for} $i=0$ {\bf up to} $rad(G),$ {\bf do} $i:$  $N_i(v_r):= \{ x \in V (G)|\ell (x) = i\}, i \in \{0, 1, 2, \ldots ,rad(G) \};$

\hskip0.2cm {\bf for} $i=0,$ then $v_r$ and $N_1(v_r)$ forming a  Fan-structure of $G.$ Let $v_r$
be the root of CCS$(G)$ and is

\hskip0.2cm marked red. {\bf For} any edge $v_{1_s}v_{1_t}$ of it, $v_{1_s},v_{1_t}\in N_1(v_r),$  {\bf if} $\mid1_s-1_t\mid\geq2,$ then
shrink the

\hskip0.2cm corresponding edge to obtain new vertex, named $v_{1v_{1_s}v_{1_t}},$ marked green and $S_{r}(v_{r})\leftarrow v_{1v_{1_s}v_{1_t}},$

\hskip0.2cm $P_{r}(v_{1v_{1_s}v_{1_t}})=v_r,$ $l_{r}(v_{1v_{1_s}v_{1_t}})=1.$

\hskip0.2cm{\bf for} $i=2$ {\bf up to} $rad(G)-1$ {\bf do}

\hskip0.4cm{\bf if} $\mid N_i(v_r)=\{v_{i_1},v_{i_2}\}\mid=2, \mid i_1-i_2 \mid\geq 2 $ and they are adjacent, then shrink the corresponding edge

\hskip0.4cm to obtain a new vertex, named $v_{iv_{i_1}v_{i_2}},$ marked green, adjacent to the new vertex $v_{1v_{1_s}v_{1_t}}$ obtained

\hskip0.4cm by the vertices in $N_{i-1}(v_r)$ which are pertinent to $N_i(v_r);$ $S_{r}(v_{1v_{1_s}v_{1_t}})\leftarrow v_{iv_{i_1}v_{i_2}},$ $P_{r}(v_{iv_{i_1}v_{i_2}})=$

\hskip0.4cm  $v_{1v_{1_s}v_{1_t}},$ $l_{r}(v_{iv_{i_1}v_{i_2}})\leftarrow i,$
$i\leftarrow i+1,$

\hskip0.4cm {\bf if} $\mid N_i(v_r)\mid\geq 3,$  then we partition them into $N_j(N_i(v)), 1\leq j \leq \mid N_i(v_r)\mid ,$ for each $j,$ the vertices of

\hskip0.4cm  $N_j(N_i(v_r))$ forming all the path edges of a Fan structure, for any edge $v_{i_s}v_{i_t}$ and $\mid i_s- i_t\mid\geq2,$ then

\hskip0.4cm  shrink the corresponding edge to obtain a new vertex, named $v_{ijv_{i_s}v_{i_t}}$ marked green, and adjacent

\hskip0.4cm  to the new vertex $v_{i-1v_{i-1_1}v_{i-1_2}}$ pertinent to $v_{i_s}v_{i_t}.$ $S_{r}(v_{1v_{1_s}v_{1_t}})\leftarrow v_{ijv_{i_s}v_{i_t}},$ $P_{r}(v_{ijv_{i_s}v_{i_t}})\leftarrow v_{1v_{1_s}v_{1_t}} $

\hskip0.4cm $l_{r}(v_{ijv_{i_s}v_{i_t}})\leftarrow i,$ $i\leftarrow i+1,$

%
%
%

\hskip0.2cm{\bf if}  $i= rad(G),$ then we partition them into $N_j(N_i(v)), 1\leq j \leq \mid N_i(v_r)\mid ,$ for each $j,$ the vertices of

\hskip0.2cm  $N_j(N_i(v_r))$ forming all the path edges of a  Fan structure, we know their common neighbor $v_s,$

\hskip0.2cm  $1\leq s\leq n$, is the central vertex of it according to {\bf the step 2}, then  $v_s$ is a vertex of CCS$(G) $ and

\hskip0.2cm marked red, which is adjacent to a vertex generated by the vertices in $N_{i-2}(v_r)$;

\hskip0.2cm If there are central vertices whose subscripts adjacent, then we shrink the corresponding edge to

 \hskip0.2cm be a green vertex connected a vertex generated by the pertinent vertices in $N_{i-2}(v_r)$; The produced

\hskip0.2cm  new vertex is adjacent to the pertinent  vertex generated by the corresponding vertices in $N_{rad(G)-2}.$

\hskip 0.2cm{\bf end}\\

The time complexity of CCS-Algorithm computing the CCS$(G)$ is $O(n^{3}),$ where $n=\mid V(G)\mid,$ for a MOP $G.$ \\

%


\section{Polynomial time algorithm for rainbow connection number of MOP}

Our main results in this section is to give a polynomial time algorithm computing rainbow coloring of MOPs with at most $3rad(G)$ colors. Specially, we can give $2rad(G)+3$ rainbow coloring for some MOPs by our algorithm.

In order to give Algorithm 3, we first give some notations for CCS$(G).$  Since CCS$(G) $ is a $v_r$ root tree and any two vertices are connected by exactly one path. Assuming that the CCS$(G)$ has $m$ leaves. Therefore there are $m$ paths between $v_r$ and it's leaves, denoted by $P_k,$ $1\leq k \leq m.$ The length $L(P_k)=:$ the number of  edges on $P_k.$ Let $P_1$ be the minimum length path of all the $m$ paths.

%

Given a MOP $G,$ algorithm 2 gives the CCS$(G).$  We perform the following operations: First,
giving the significant Fan-structures, whose central vertices are corresponding to the vertices of
CCS$(G),$ of $G,$ pertinent to  CCS$(G):$ (1) we choose $v_r$ and it's neighbors as the first
 Fan-structure, denoted by $Fan_r.$ If a vertex of CCS$(G)$ is red, then it is a vertex of $G$ and
  forms up a Fan-structure as a central vertex with it's neighbors;
  If a vertex $v_{iv_{i_1}v_{i_2}}$ is green,
  then we choose the vertex corresponding to the $\max\{\mid v_{i_1}\mid,\mid v_{i_2}\mid\}$
  as central vertices of $G,$ which form Fan-structures with it's neighbors.
  If there are neighbors in
   $N(\max\{degree( v_{i_1}),degree( v_{i_2})\})\backslash N(\min\{degree( v_{i_1}),degree( v_{i_2})\}), $
   then we choose $v_{i_1},v_{i_2}$ as the central vertices of Fan-structures in $G.$ Second,
   giving a rainbow coloring using at most $3rad(G)$ colors.\\


 \noindent{\bf Algorithm 3 giving rainbow coloring of MOP:}\\

\noindent{\bf Input:} A  MOP $G.$\\
\noindent{\bf Output:} A rainbow coloring of $G.$

{\bf Step 1:} Giving the CCS$(G)$ of $G$ by {\bf Algorithm 2}.

{\bf Step 2:} Give three color sets: $C_1=\{ 7,8,\cdots rad(G)+4\}$ and $C_2=\{ rad(G)+5,rad(G)+7,\cdots$

\hskip0.2cm $ 2rad(G)+2\},$ $C_3=\{2rad(G)+3,\cdots 3rad(G)\}.$

{\bf Step 3:} We produce the corresponding two path for every $P_k,1\leq k\leq m$
by the procedures of


\hskip0.2cm {\bf Theorem 3.4.} Now we color the edges of the shorter paths by $C_1$ and the other by the the unused

\hskip0.2cm  minimum colors in $C_2$ and $C_3.$

%

%
%
%

{\bf Step 4:} In this step, we give a $2rad(G)+2+c,$ where $c\leq rad(G)-2,$ rainbow coloring of $G.$

\hskip0.2cm (1) Choose $r$ and it's neighbors as the first Fan, denoted by $Fan_r,$  Whose spoke edges are colored

\hskip0.2cm by colors $4,5$ alternatively and according to the clockwise around the central vertex and uncolored

\hskip0.2cm path edges are colored by $6.$

\hskip0.2cm (2){\bf begin}

\hskip0.2cm {\bf for} $k=1$ {\bf up to} $m$ {\bf do}

\hskip0.2cm{\bf if} $v_e\in V(P_k)$ and is a red vertex of CCS$(G),$ then it is a vertex of $G$ and forms up a Fan-structure

\hskip0.2cm as a central vertex with it's neighbors. According to {\bf the step 2 of Algorithm 2}, we know the

\hskip0.2cm Fan-structure. Whose uncolored spoke edges are colored by colors $1,2$ alternatively and according

\hskip0.2cm to the clockwise around the central vertex and uncolored path edges are colored by $3;$

\hskip0.2cm{\bf if} $v_e\in V(P_k)$ and is a green vertex $v_{{ij}v_{i_1}v_{i_2}},$ {\bf first} choose the vertex corresponding to the

\hskip0.2cm $\max\{degree( v_{i_1}),degree( v_{i_2})\}$ and it's neighbors which forming a Fan-structure in $G,$ {\bf second}, if

\hskip0.2cm there are neighbors in $N(\min\{degree( v_{i_1}),degree( v_{i_2})\})\backslash N(\max\{degree( v_{i_1}),degree( v_{i_2})\}) $ then we

\hskip0.2cm choose $v_{i_1},v_{i_2}$ as the central vertices of Fan-structures in $G.$  According to {\bf the step 2 of Algorithm

\hskip0.2cm 2}, we know the Fan-structures. They are colored with $1,2,3$ by the method as above shown.

\hskip 0.2cm{\bf end}\\

\noindent{\bf Theorem 4.1.} If $G$ is a MOP, then the edge coloring given by algorithm 3 is a rainbow coloring of $G.$ Moreover it uses colors at most $3rad(G).$\\

\noindent{\bf Proof.}

{\bf Case 1:} For any two vertices, which locate in a Fan-structure whose central vertex corresponding the vertex of CCS$(G).$ Because it's spoke edges are colored by colors $1,2$ alternatively and according to the clockwise around the central vertex and path edges are colored by $3,$ other colors, if exist, belong to $C_1,C_2$ or $C_3.$ Obviously, there is a rainbow path between them.

If two vertices locate in different Fan-structures whose central vertices corresponding two vertices, which are situated on a $v_r-path,$  of CCS$(G).$ Then there is a rainbow path, which use colors of $C_1,C_2$ or $C_3$ connecting their central vertices in $G.$ Moreover the two Fan-structures are colored with the above method, therefore the two vertices can be rainbow connected under the coloring given by Algorithm 3.

 {\bf Case 2:} We know that $v_r$ with it's neighbors forming a Fan-structure, whose spoke edges are colored by $4,5$ alternatively and according to the clockwise around it and uncolored path edges are colored by $6.$ Obviously any two vertices of the Fan structure are rainbow connected under the coloring.


 Any two vertices whose central vertices are saturated on different $v_r-paths$ are rainbow connected through $Fan_{v_r}.$  Since there are two edge disjoint rainbow paths connecting the central vertices of the corresponding Fan-structures and the Fan $Fan_{v_r}.$

For  any longer path of Step 3, since their first edges are colored by 4 or 5 and other edges no more than $2(rad(G)-2),$ which can be rainbow colored by the colors of $C_2$ and $C_3.$   $\hfill \Box$\\


\section{Examples of rainbow connection algorithms for MOPs}

X.Deng  et al. $\cite{XD}$ gave a polynomial time algorithm to compute the rainbow connection number of MOPs, but only obtain an compact upper bound, by the Maximal fan partition method. By the method, they proved that $rc(Figure 13)\leq 6.$ \\

 \begin{figure}[ht]
\begin{center}
\includegraphics[width=5cm,totalheight=30mm]{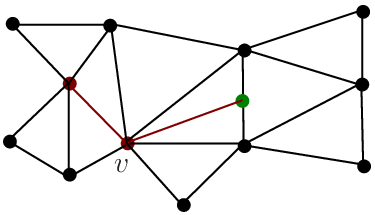}
\end{center}
\caption{$\mathrm{An~ example}$}\label{fig11}
 \end{figure}

Obviously, we know that the vertex $v$ is the central vertex of Figure \ref{fig11} and the red and green vertices and red edges are CCS$(Figure \ref{fig11}),$ then we know that $rc(Figure \ref{fig11})\leq rc(\mbox{CCS}(Figure \ref{fig11}))+3=5$ by the coloring of {\bf Algorithm 3}.

 Figure \ref{fig8} is a MOP with diameter $3,$ where the Maximal fan structures having $9$ vertices and the red vertices and edges give the CCS$(Figure\ref{fig8}).$ When the Fan structures are replaced by any $F_n$ with $n\geq 9,$ the resulting graphs can be rainbow colored by $4$
colors with the same method as Figure \ref{fig8} showing. If we color the graph by the Maximal fan partition method, then the colors needed increase to 9 with the rising of $n.$\\

\noindent{\bf Definition $5.1^{\cite{CDRV}}.$  ({\bf Two-way dominating set}).} A dominating set $D$ in a graph $G$ is called
a two-way dominating set, if every pendant vertex of $G$ is included in $D.$ In addition,
if $G[D]$ is connected, we call $D$ a connected two-way dominating set.\\

\noindent{\bf Theorem $5.2^{\cite{CDRV}}.$ } If $D$ is a connected two-way dominating set in a graph $G,$ then $rc(G)\leq rc(G[D])+3.$\\

By the two-way dominating set and induction on radius, they proved the following result:\\

\noindent{\bf Theorem $5.3^{\cite{CDRV}}.$} If $G$ is a bridge-less chordal graph, then $rc(G)\leq 3¡¤rad(G).$ Moreover,
there exists a bridge-less chordal graph with $rc(G)=3¡¤rad(G).$\\

For Figure \ref{fig11}, the minimum connected two-way dominating set is a $4-length$ path, thus
  $rc(Figure \ref{fig11})\leq 3rad(Figure \ref{fig11})=6$ by the theorem 5.3.

Another example Figure \ref{Fan70}, we know that Figure \ref{Fan70} is a MOP, obviously a chordal graph, and has a hamiltonian cycle with $40$ vertices. By Theorem 1.1, we know that  $rc(Figure\ref{Fan70})\leq 12.$ It is smaller than it's $3rad(Figure\ref{Fan70})=15$ and bigger 2 than it's diameter.

If a MOP $G$ is $Fan_{n},n\geq 7,$ then $rad(G)=1,rc(G)=3rad(G)=3.$ Therefor the upper bound of the algorithm given is sharp for MOPs, so is the {\bf Teorem 5.3.} But above all, the bound gived by our algorithm is better to the theorem 5.3 obtained for some MOPs. For example, Figure \ref{Fan70} showing. \\

\section{Concluding remarks}
 Recall algorithm 1, we know that the MOPs $Lad_d$ and $Lad_{d^+},$ with diameter $d,$ have rainbow connection number $d.$ Producing the graphs and giving their strong rainbow connection numbers take time at most $O\left(d\right).$

 In algorithm 2: The Hamiltonian cycle $C_G$ of MOP $G$ can be obtained by a
linear time algorithm presented in $\cite{SM}$ through the canonical
representation of $G.$ Then select any vertex of $G$ as the initial
vertex of $C_G,$ we can obtain a Hamiltonian degree sequence
in time $O\left(n^3\right),$ where $n=n(G).$ Note the time of Step 2
is $O(n+e),$ where $e=e(G)$ and  $e(G)=2n-3,$ which is $O(n),$ Step 3 takes time at most $O\left(n^3\right) $ and Step 4 at most $O\left(n\right).$ For MOP, Step 5 takes time $O(n+e),$ which is $O(n).$ Step 6 at most $O(n^2).$

For algorithm 3: The time of Step 3 is no more than $O(n),$   Step 4 uses time at most $O(n^2).$ If $G$ is a MOP, the algorithm gives a tight upper bound of $rc(G)$ which is no more than $2rad(G)+2+c,$ where $0\leq c\leq rad(G)-2$ and needs time at most $O\left(n^3\right).$

In the future, we will consider the rainbow connection numbers for outerplanar and general planar graphs on algorithm aspect.
It is interesting to study the rainbow connection for planar graphs on algorithm aspect.\\

\end{document}